\documentclass{article}
\usepackage{amssymb}

\newtheorem{theorem}{Theorem}
\newtheorem{proposition}[theorem]{Proposition}
\newtheorem{corollary}[theorem]{Corollary}

\newtheorem{remark}[theorem]{Remark}
\newcommand{\dem}{\noindent{\bf Proof. }}

\title{{ Dynamics of a rational system of difference equations in the plane }}

\author{Ignacio Bajo, Daniel
Franco and Juan Per\'an}
\begin{document}

\maketitle

\begin{center} {\bf Abstract}\end{center}

\noindent We consider a rational system of first order difference equations in the plane with four  parameters such that all fractions have a common denominator. We study, for the different values of the parameters, the global and local properties of the system. In particular, we discuss  the boundedness and the asymptotic behavior of the solutions, the existence of periodic solutions and the stability of equilibria.

\bigskip

\noindent {\it Keywords:}
Rational system in the plane; bounded solutions; periodic solutions; local stability; global stability.

\noindent {\it AMS Class.:} 39A22;39A23;39A30

\section{Introduction}
In  recent years, rational difference equations have attracted the attention of many researchers for varied reasons. On the one hand, they provide examples of non-linear equations which are, in some cases, treatable but whose dynamics present some new features with respect to the linear case. On the other hand, rational equations frequently appear in some biological models and, hence, their study is of interest also due to their applications. A good example of both facts are Ricatti difference equations; the richness of the dynamics of Ricatti equations is very well-known (see, for instance, \cite{cull} and \cite{elaydi}) and a particular case of these equations provides the classical  Beverton-Holt model on the dynamics of exploited fish populations \cite{BH}. Obviously, higher order rational difference equations and  systems of rational equations have also been widely studied but still have many aspects to be investigated. The reader can find in the following books \cite{1},\cite{2},\cite{3}, and the works cited therein many results, applications and open problems on higher order equations and rational systems.

A preliminar study of planar rational systems in the large can be found in the paper \cite{ladas09} by Camouzis et al.  In such work, they give some results and provide some open questions for systems of equations of the type
$$\left.\begin{array}{l}
x_{n+1}=\displaystyle\frac{\alpha_1+\beta_1x_n+\gamma _1y_n}{A_1+B_1x_n+C_1 y_n}\vspace{0.2cm}\\
y_{n+1}=\displaystyle\frac{\alpha_2+\beta_2x_n+\gamma _2y_n}{A_2+B_2x_n+C_2 y_n}\end{array}\right\}\, ,\quad n=0,1,\dots
$$
where the parameters are taken to be non-negative. As shown in the cited paper, some of those systems can be reduced to some Ricatti equations or to some previously studied second-order rational equations. Further, since for some choices of the parameters one obtains a system which is equivalent to the case with some other parameters, Camouzis et al. arrived at a
 list of 325 non-equivalent systems to which it should be focused the attention. They list such systems as pairs $(k,l)$ where $k$ and $l$ make reference to the number of the corresponding equation in their Tables 3 and 4.

In this paper,  we deal with the rational system labelled as (21,23) in \cite{ladas09}. Note that for non-negative coefficients such system is neither cooperative nor competitive but it has the particularity that denominators in both equations are equal. This allows us to use some of the techniques developed in \cite{b-l} to completely obtain the solutions and give a nice description of the dynamics of the system. In principle, we will not restrict ourselves to the case of non-negative parameters, although this case will be considered in detail in the last section. Hence, we will study the general case of the system
\begin{equation}\label{sistema}\left.\begin{array}{l}
x_{n+1}=\displaystyle\frac{\alpha_1+\beta_1x_n}{y_n}\vspace{0.2cm}\\
y_{n+1}=\displaystyle\frac{\alpha_2+\beta_2x_n}{y_n}\end{array}\right\}\, ,\quad n=0,1,\dots ,
\end{equation}
where the parameters $\alpha_1,\alpha_2,\beta_1,\beta_2$ are given real numbers and the initial condition $(x_0,y_0)$ is
an arbitrary vector of $\mathbb{R}^2$. It should be noticed that when $\alpha_1\beta_2=\alpha_2\beta_1$ the system can be reduced to a Ricatti equation (or it does not admit any complete solution, which occurs for $\alpha_2=\beta_2=0$) and therefore these cases will be neglected. Since we will not assume non-negativeness for neither the coefficients nor the initial conditions, a forbidden set will appear. We will give an explicit characterization of the forbidden set in each case. Obviously,  all the results concerning solutions that we will state in the paper are to be apply only to complete orbits.
We will focus our attention in three aspects of the dynamics of the system: the boundedness character and asymptotic behavior of its solutions, the existence of periodic orbits (and, in particular, of prime period-two solutions) and the stability of the equilibrium points. It should be remarked that, depending on the parameters, they may appear asymptotically stable fixed points, stable but not asymptotically stable fixed points, non-attracting unstable fixed points and attracting unstable fixed points.

The paper is organized, besides this introduction, in three sections. Section 2 is devoted to some preliminaries and some results which can be mainly deduced from the general situation studied in \cite{b-l}. Next, we study the case $\beta_2=0$ since such assumption yields the uncoupled globally 2-periodic equation $y_{n+1}=\alpha_2/y_{n}$ and the system is reduced to a linear first order equation with 2-periodic coefficients; this will be our section 3 below. The main section of the paper is section 4, where we give the solutions to the system and the description of the dynamics in the general case $\beta_2\neq 0$. We finish the paper by describing the dynamics in the particular case where the coefficients and the initial conditions are taken to be non-negative.

\section{Preliminaries and first results}
Systems of linear fractional difference equations $X_{n+1}=F(X_n)$ in which denominators are common for all the components of $F$ have been studied in \cite{b-l}. If one denotes by $q$ the mapping given by $q(a_1,a_2,\dots,a_{k+1})=(a_1/a_{k+1},a_2/a_{k+1},\dots, a_k/a_{k+1})$ for $(a_1,a_2,\dots,a_{k+1})\in {\mathbb R}^{k+1}$ with $a_{k+1}\neq 0$ and $\ell: {\mathbb R}^k\to {\mathbb R}^{k+1}$ is given by $\ell ( a_1,a_2,\dots,a_{k})=(a_1,a_2,\dots,a_{k},1)$, it is shown in such work that the system can be written in the  form $X_{n+1}=q\circ A\circ\ell (X_n)$, where $A$ is a $(k+1)\times (k+1)$ square matrix constructed with the coefficients of the system. In the special case of our system (\ref{sistema}) one actually has
$$\left(\begin{array}{c} x_{n+1}\\y_{n+1}\end{array}\right)=q\circ
\left(\begin{array}{ccc}\beta_1&0&\alpha_1\\\beta_2&0&\alpha_2\\0&1&0\end{array}\right)\left(\begin{array}{c} x_{n}\\y_{n}\\1\end{array}\right).$$
This  form of the system let us completely determine its solutions in terms of the powers of the associated matrix
\begin{equation}\label{matriz}
A=\left(\begin{array}{ccc}\beta_1&0&\alpha_1\\\beta_2&0&\alpha_2\\0&1&0\end{array}\right).
\end{equation}
Actually, the explicit solution to the system with initial condition $(x_0,y_0)$ is given by
\begin{equation}\label{expmat}(x_{n+1},y_{n+1})^t=q\circ A^n (x_0,y_0,1)^t,\end{equation}
where $M^t$ stands for the transposed of a matrix $M$. Therefore, our system can be completely solved and the solution starting at $(x_0,y_0)$ is just the projection by $q$ of the solution of the linear system $X_{n+1}=AX_n$ with initial condition $X_0=(x_0,y_0,1)^t$ whenever such projection exists.

\begin{remark} When such projection does not exist, then $(x_0,y_0)$ lies in the forbidden set. Clearly, this may only happen when for some $n\geq 1$ one has
$$(0,0,1)A^n(x_0,y_0,1)^t=0.$$
Therefore, if $a_i(n)\in {\mathbb R}$, $0\leq i\leq 2$ are such that
$A^n=a_0(n)I+a_1(n)A+a_2(n)A^2$, then one obtains immediately that
the forbidden set is given by  the following union of lines
$$\mathbf{F}=\bigcup_{n\geq 1}\{\ (x_0,y_0)\in {\mathbb R}^2\, :\, a_1(n)y_0+a_2(n)\beta_2x_0+a_2(n)\alpha_2+a_0(n)=0\  \}.$$

The explicit calculation of $a_i(n),0\leq i\leq 2$ for each $n\geq 3$ may be done in several ways. For instance, one has that $a_0(n)+a_1(n)x+a_2(n)x^2$ is the remainder of the division of $x^n$ by the characteristic polynomial of $A$. Further, by elementary techniques of Linear Algebra one can also compute them in terms of the eigenvalues of $A$ (an approach using the solutions to an associated linear difference equation may be seen in \cite{elaydi2}). \end{remark}

\begin{remark} As mentioned in the introduction, all through the paper we will consider that \begin{equation}
\label{detno0}
\beta_2 \alpha_1 \ne \beta_1 \alpha_2.
\end{equation}
(this is to say, that the matrix $A$ is non-singular) since the cases with $\beta_2 \alpha_1 = \beta_1 \alpha_2$ may be reduced to a single Ricatti equation. Actually, if $\alpha_2=\beta_2=0$, then the system does not admit any complete solution, whereas, for $\alpha_2\ne 0$ or $\beta_2\ne 0$, one has that there exists a constant $C$ such that
$\alpha_1=C\alpha_2$ and $\beta_1=C\beta_2$ and hence the first equation of the system may be substituted by $x_{n+1}=Cy_{n+1}$ and then the second one  reduces to the Ricatti equation
$$y_{n+2}=\frac{\alpha_2+\beta_2Cy_{n+1}}{y_{n+1}}\, ,\,\, n=0,1,\dots$$ with initial condition $y_1=\frac{\alpha_2+\beta_2x_0}{y_0}.$
\end{remark}

\bigskip

Our main goal will be to give a description of the dynamics of the system in terms of the eigenvalues of the associated matrix $A$ given in (\ref{matriz}). We begin with the following result concerning 2-periodic solutions which is the particularization to our system of the analogous general result given in Theorem 3.1 and Remark 3.1 of \cite{b-l}.

\begin{proposition} \label{sol2p} Let us consider the system (\ref{sistema}) with $\alpha_1\beta_2\neq\alpha_2\beta_1$. One has:
\begin{enumerate}
\item If $\beta_2\neq 0$, then there are exactly as many equilibria as distinct real eigenvalues of the matrix $A$. More concretely, for each real eigenvalue $\lambda$ one gets the equilibrium $\left(\frac{\lambda^2-\alpha_2}{\beta_2} , \lambda\right).$
\item When $\beta_2=0$, one has:
\begin{enumerate}
\item if $\alpha_2<0$, then there are no fixed points;
\item if $0<\alpha_2 \neq\beta_1^2$, then there are two fixed points at $\left(\frac{\alpha_1}{\sqrt{\alpha_2}-\beta_1},\sqrt{\alpha_2}\right)$ and
$\left(\frac{-\alpha_1}{\sqrt{\alpha_2}+\beta_1},-\sqrt{\alpha_2}\right)$;
\item  if $\alpha_2
=\beta_1^2$ and $\alpha_1\neq 0$, then the only equilibrium point is $\left(\frac{-\alpha_1}{2\beta_1}, -\beta_1\right)$;
\item  if $\alpha_2
=\beta_1^2$ and $\alpha_1= 0$, then there is an isolated fixed point $(0,-\beta_1)$ and a whole line of equilibria $(x_0,\beta_1)$.
\end{enumerate}
\item There exist periodic solutions of prime period 2 if and only if $\alpha_1\beta_2=0$.
\end{enumerate}
\end{proposition}
\dem  As stated in \cite{b-l}, a point $(a,b)\in {\mathbb R}^2$ is an equilibrium if and only if $(a,b,1)$ is an eigenvector of the associated matrix $A$. When $\beta_2\neq 0$, it is straightforward to prove that for each real eigenvalue $\lambda$ the vector $\left(\frac{\lambda^2-\alpha_2}{\beta_2} , \lambda , 1\right)$ is an eigenvector. In the case $\beta_2=0$, the equilibrium points can be easily computed directly from the equations
$\alpha_2=y^2$, $\alpha_1+\beta_1 x= x y.$

For the proof of affirmation (3) it suffices to bear in mind that, according to \cite{b-l}, the existence of prime period-two solutions is only possible when the associated matrix $A$ has an eigenvalue $\lambda$ such that $-\lambda$ is also an eigenvalue. Since $A$ is a $3\times 3$ square matrix, this obviously implies that the trace of $A$ is also an eigenvalue. Hence, $\beta_1$ is an eigenvalue but this is only possible if $\alpha_1\beta_2=0.$ If $\alpha_1=0$, then the initial condition $(0,y_0)$ gives a prime period 2 solution whenever $y_0^2\neq \alpha_2$, whereas, if $\alpha_1\neq 0$ and $\beta_2= 0$, a direct calculation shows that the solution with initial conditions $(0,-\beta_1)$ is periodic of prime period 2.  $\square$

We now study the stability of fixed points in some of the cases. Recall that a fixed point of our system $(x^*,y^*)$ always verifies $y^*=\lambda$ for some real eigenvalue $\lambda$ of the matriz $A$. We will say in such case that the fixed point  $(x^*,y^*)$ is {\it associated} to $\lambda$.

\begin{proposition} \label{esta-linear} Let us consider the system (\ref{sistema}) with $\alpha_1\beta_2\neq\alpha_2\beta_1$. Let $\rho (A)$ be the spectral radius of the matrix $A$ given in (\ref{matriz}) and let $\lambda$ be an eigenvalue of $A$.
\begin{enumerate}
\item If $|\lambda|<\rho (A)$, then the associated equilibrium is unstable.
\item If $|\lambda|=\rho (A)$ and all the eigenvalues of $A$ whose modulus is $\rho (A)$ are simple, then the associated fixed point is stable. Further, if in this case $\lambda$ is the unique eigenvalue whose modulus is $\rho (A)$, then it is assymptotically stable.
\end{enumerate}
\end{proposition}
\dem  The Jacobian matrix of the map $F(x,y)=\left(\frac{\alpha_1+\beta_1 x}{y},\frac{\alpha_2+\beta_2 x}{y}\right)$ at a fixed point $(x^*,y^*)$ is given by
$$DF(x^*,y^*)=\left(\begin{array}{ccl} \beta_1/y^*& \, &-x^*/y^*\\ \beta_2/y^*& &-1\end{array}\right).$$
Consider an eigenvalue $\lambda$ of $A$ and let $\lambda_2,\lambda_3$ be the other (non-necessarily different) eigenvalues of $A$. Let us show that the eigenvalues of the Jacobian matrix at a fixed point associated to $\lambda$ are just $\frac{\lambda_2}{\lambda}$ and $\frac{\lambda_3}{\lambda}$. The result is trivial when $\beta_2=0$ since the eigenvalues of $A$ are $\beta_1$ and $\pm\sqrt{\alpha_2}$ and fixed points are always associated to one of the eigenvalues $\pm\sqrt{\alpha_2}$. If $\beta_2\neq 0$, then $x^*=\frac{{\lambda}^2-\alpha_2}{\beta_2}$ and $y^*=\lambda$ and, therefore, one obtains:
\begin{eqnarray*}
\mbox{trace}(DF(x^*,y^*))&=& \frac{\beta_1-\lambda}{\lambda}=\frac{\lambda_2+\lambda_3}{\lambda}\\
\mbox{det}(DF(x^*,y^*))& = & \frac{-\beta_1\lambda+\lambda^2-\alpha_2}{\lambda^2}=\frac{\mbox{det}(A)}{\lambda^3}=\frac{\lambda_2\lambda_3}{\lambda^2},
\end{eqnarray*}
showing that the eigenvalues of $DF(x^*,y^*)$ are as claimed. Now, the first statement follows at once since, if $|\lambda|<\rho (A)$, then at least one of the eigenvalues of $DF(x^*,y^*)$ lies outside the unit circle. Moreover, when $|\lambda|=\rho (A)$ and it is the unique eigenvalue with such property, then the eigenvalues of $DF(x^*,y^*)$ are inside the (open) unit ball and hence the equilibrium $(x^*,y^*)$ is assymptotically stable, which proves the second part of (2).

For the proof of the first part of (2) let us recall that if $(x^*,y^*)$ is a fixed point of (\ref{sistema}) associated to the real eigenvalue $\lambda$, then $X^*=(x^*,y^*,1)^t$ is a fixed point of the linear system $X_{n+1}=\frac{1}{\lambda}A X_n$. The eigenvalues of the matrix $M=\frac{1}{\lambda}A $ are obviously $1,\frac{\lambda_2}{\lambda}$ and $\frac{\lambda_3}{\lambda}$. Since the eigenvalues of $A$ having modulus  $\rho (A)$ are simple, so are the eigenvalues of $M$ having modulus $1$. Therefore, the fixed point $X^*$ is stable \cite[Th. 4.13]{elaydi}. Now, the stability of $(x^*,y^*)$ follows at once from (\ref{expmat}) and the continuity of $q$ in the semi-space $z>0$. 
$\square$

\section{Case $\beta_2=0$.}

Recall that, since we are assuming that inequality (\ref{detno0}) holds, we have  $\beta_1 \alpha_2 \ne 0$. In this case, the forbidden set of the system reduces to the line $y=0$. Since $\beta_2=0$ the second equation of the system becomes the uncoupled equation
$$y_{n+1}=\frac{\alpha_2}{y_n},$$ which, as far as $\alpha_2\ne
0$,  for each initial condition $y_0\ne 0$ gives
\begin{equation}\label{solyb20} y_n=
\left \{
\begin{array}{ll}
\medskip
y_0 & \mbox{ for even $n$},\\
\displaystyle\frac{\alpha_2}{y_0}  & \mbox{ for odd $n$.}
\end{array}
\right .
\end{equation}
 Substituting such values in the first equation
of the system we obtain a first order linear difference equation with 2-periodic coefficients whose solution is given by $x_1=(\alpha_1+\beta_1x_0)/y_0$ and, for $n>1$,
\begin{equation}
\label{solxb20}
x_n=
\left \{
\begin{array}{ll}\medskip
\displaystyle \left ( \frac{{\beta_1}^2}{\alpha_2}\right )^{\frac{n}{2}} \left [ x_0 + \displaystyle \frac{\alpha_1
(\beta_1 + y_0)}{\alpha_2} \sum_{k=1}^{\frac{n}{2}} \left ( \frac{\alpha_2}{{\beta_1}^2}\right )^k \right ] &  \mbox{for
even $n$},\\
\displaystyle \frac{ \alpha_1}{y_0 }+ \frac{\beta_1}{y_0} \left( \frac{{\beta_1}^2}{\alpha_2}\right )^{\frac{n-1}{2}}
\left [ x_0  +  \frac{\alpha_1 (\beta_1 + y_0)}{\alpha_2}\sum_{k=1}^{\frac{n-1}{2}}
\left(\frac{\alpha_2}{{\beta_1}^2}\right)^{k}  \right ]  &  \mbox{for odd $n$.}
\end{array}
\right .
\end{equation}

Hence, we have proved the following.

\begin{proposition}
\label{tb20}
If $\beta_2=0$ and $\beta_2\alpha_1\ne \beta_1 \alpha_2$, then the system (\ref{sistema}) is solvable for any initial condition $(x_0,y_0)$
with $y_0\ne 0$ and the solution $(x_n,y_n)$ is given by (\ref{solyb20}) and (\ref{solxb20}) where, explicitly, one has:\\

(1) If  $\alpha_2=\beta_1^2$, then for $n>1$
\[ x_n=
\left \{
\begin{array}{ll}\medskip
\displaystyle  x_0 - \frac{\alpha_1 (\beta_1 + y_0)n}{2{\beta_1}^2} &  \mbox{ for even $n$},\\
\displaystyle \frac{ \alpha_1 }{y_0}+\frac{ \beta_1x_0}{y_0}-  \displaystyle \frac{\alpha_1 (\beta_1 + y_0)(n-1)}{2\beta_1y_0}    &
\mbox{ for odd $n$}.
\end{array}
\right .
\]

(2) If  $\alpha_2\neq\beta_1^2$, then for $n>1$
$$
x_n=
\left \{
\begin{array}{ll}\medskip
\displaystyle \left ( \frac{{\beta_1}^2}{\alpha_2}\right )^{\frac{n}{2}} \left [ x_0 + \displaystyle \frac{\alpha_1
(\beta_1 + y_0)}{{\beta_1}^2-\alpha_2} \left ( 1-\left ( \frac{\alpha_2}{{\beta_1}^2}\right )^{\frac{n}{2}} \right )
\right ] &  \mbox{for even $n$},\\
\displaystyle \frac{ \alpha_1}{y_0 }+ \frac{\beta_1}{y_0} \left ( \frac{{\beta_1}^2}{\alpha_2}\right )^{\frac{n-1}{2}}
\left [ x_0 + \displaystyle \frac{\alpha_1 (\beta_1 + y_0)}{{\beta_1}^2-\alpha_2} \left ( 1-\left (
\frac{\alpha_2}{{\beta_1}^2}\right )^{\frac{n-1}{2}} \right ) \right ] &  \mbox{for odd $n$.}
\end{array}
\right .
$$
\end{proposition}

\smallskip

From the proposition above one can easily derive the following result which completely describes the asymptotic behaviour of the solutions to the system.

\smallskip

\begin{corollary} Let us consider $\beta_2=0$ and $\beta_1\alpha_2\neq 0$.
\begin{enumerate}
\item When ${\beta_1}^2 = \alpha_2$ one has:
\begin{enumerate}
\item If $\alpha_1\ne 0$, then every solution to the system is unbounded except those with initial condition $(x_0,-\beta_1)$, which are 2-periodic.
\item If $\alpha_1= 0$, the system is globally 2-periodic.
\end{enumerate}
\item If ${\beta_1}^2 = -\alpha_2$, then the system (\ref{sistema}) is globally 4-periodic. Further, the solution corresponding with the initial condition $(x_0,y_0)$ is of prime period 2 if and only if  $2\beta_1^2x_0+\alpha_1(\beta_1+y_0)=0$.
\item If ${\beta_1}^2 \ne |\alpha_2|$, then the solutions with initial condition
 $\left(\displaystyle \frac{\alpha_1(\beta_1+y_0)}{\alpha_2-\beta_1^2},y_0\right)$ are period-two solutions. Moreover,
\begin{enumerate}
\item If ${\beta_1}^2 > |\alpha_2|$, then any other solution to the system (\ref{sistema}) is unbounded.
\item If ${\beta_1}^2 < |\alpha_2|$, then any other solution of (\ref{sistema}) is bounded and tends to one of the
period-two solutions described above.
\end{enumerate}\end{enumerate}
\end{corollary}
\dem 
The proof is a straightforward consequence of the explicit formulas for $x_n$ and $y_n$ given in Proposition \ref{tb20}. It should, however, be mentioned that the globally periodicity of the system in the case ${\beta_1}^2 = -\alpha_2$ can be easily seen since the associated matrix $A$ given by (\ref{matriz}) in such case  verifies  $A^4=\beta_1^4 I$, where $I$ stands for the identity matrix. Actually, a simple calculation proves that the solution starting at $(x_0,y_0)$ is the 4-cycle
$$ \left\{(x_0,y_0),\left(\frac{\alpha_1+\beta_1 x_0}{y_0},\frac{-\beta_1^2}{y_0}\right),\left( -x_0-\frac{\alpha_1(\beta_1+y_0)}{\beta_1^2},y_0\right),\left(\frac{-\beta_1^2x_0+\alpha_1y_0}{\beta_1y_0}, \frac{-\beta_1^2}{y_0}\right)\right\},$$
which is obviously 2-periodic if and only if $x_0=-x_0-\frac{\alpha_1(\beta_1+y_0)}{\beta_1^2}$.$\square$

From the above result and Proposition \ref{esta-linear} one easily gets the following information about the stability of the fixed points.
\begin{corollary}
Let us consider $\beta_2=0$ and $\beta_1\alpha_2\neq 0$.
\begin{enumerate}
\item If ${\beta_1}^2 = \alpha_2$, then
\begin{enumerate}
\item for $\alpha_1\ne 0$ the unique fixed point of (\ref{sistema}) is unstable;
\item for $\alpha_1= 0$ every fixed point of (\ref{sistema}) is stable but not asymptotically stable.
\end{enumerate}
\item If ${\beta_1}^2 \ne \alpha_2 > 0$, then
\begin{enumerate}
\item for ${\beta_1}^2 > \alpha_2$ both fixed points of (\ref{sistema}) are unstable;
\item for ${\beta_1}^2 < \alpha_2$ the fixed points of (\ref{sistema}) are stable but not asymptotically stable.
\end{enumerate}\end{enumerate}
\end{corollary}

\section{Case $\beta_2\ne 0$.}

\begin{proposition}
\label{solb2no0}
 Suppose $\beta_2\ne 0$ and  $(x_0,y_0)$ is an initial condition not belonging to the forbidden set $\mathbf{F}$. In such
 case the solution of system (\ref{sistema}) is given by
\begin{equation}\nonumber
x_n= \displaystyle \frac{v_{n+1}}{v_{n-1}}\displaystyle \frac{1}{\beta_2}-\displaystyle \frac{\alpha_2}{\beta_2}\, ,\qquad y_{n}=\displaystyle \frac{v_{n}}{v_{n-1}},
\end{equation}
where $v_n$ is the unique solution of the linear difference equation
\begin{equation}
\label{ecuacionlineal}
v_{n+3}-\beta_{1}v_{n+2}-\alpha_{2}v_{n+1}+\left(
\beta_{1}\alpha_{2}-\beta_{2}\alpha_{1}\right) v_{n}=0
\end{equation}
with initial conditions $v_{-1}=1$, $v_0=y_0$ and $v_1=\beta_2 x_0+\alpha_2$.
\end{proposition}
\dem  As we have seen in section 2, the solution to system (\ref{sistema}) starting at a point $(x_0,y_0)$ not belonging to the forbidden set is just the projection by $q$ of the solution of the linear system
$(u_{n+1},v_{n+1},w_{n+1})^t=A(u_n,v_n,w_n)^t$ with initial condition $(x_0,y_0,1)^t$ where $A$ is given by (\ref{matriz}).
Since the third equation of such linear systems reads $w_{n+1}=v_{n}$, it can be reduced to the planar linear system of second order equations
\begin{equation}
\label{sistemaenuv}
\begin{array}{ll}
u_{n+1} & =\beta_{1}u_{n}+\alpha_{1}v_{n-1}, \\
v_{n+1} & =\beta_{2}u_{n}+\alpha_{2}v_{n-1},
\end{array}
\end{equation}
and hence, if $(u_n,v_n)$ is the solution to (\ref{sistemaenuv}) obtained for the initial conditions $(u_0,v_0,v_{-1})=(x_0,y_0,1)$, then the solution of our rational system for the initial values $(x_0,y_0)$ will be
\begin{equation}\label{solufinal}
\begin{array}{ll}
x_{n+1}=u_n/v_{n-1} \, ,\quad y_{n+1}=v_n/v_{n-1}.
\end{array}
\end{equation}
It is clear that for $\beta_2\neq 0$ we have that $u_n$ can be completely determined  by (\ref{sistemaenuv}) in terms of $v_{n+1}$ and $v_{n-1}$ and hence it suffices to solve the third order linear equation
\begin{equation}\nonumber
v_{n+3}-\beta_{1}v_{n+2}-\alpha_{2}v_{n+1}+\left(
\beta_{1}\alpha_{2}-\beta_{2}\alpha_{1}\right) v_{n}=0
\end{equation}
trivially deduced from (\ref{sistemaenuv}) and substitute the corresponding values in (\ref{solufinal}) to obtain the result claimed.$\square$

In the following results we shall discuss the behavior of the solutions to (\ref{sistema}) by using Proposition \ref{solb2no0}. We shall consider three different cases depending on the roots of the characteristic
polinomial of the linear equation (\ref{ecuacionlineal}). Recall that such roots are also the (possibly complex) eigenvalues of the matrix $A$ given in (\ref{matriz}).

 From Proposition \ref{solb2no0} we see that the asymptotic behavior of the solutions of system (\ref{sistema}) will depend on the asymptotic behavior of the sequences $\frac{v_n}{v_{n-1}}$,
 being $v_n$ solutions of the linear difference equation  (\ref{ecuacionlineal}). The Theorem of Poincar\'e (Theorem 8.9 in \cite{elaydi}) establishes a general result for the existence of $\lim_{n\to \infty} \frac{v_n}{v_{n-1}}$. In our case, since equation (\ref{ecuacionlineal}) has constant coefficients, we can directly do the calculations, even in the cases not covered by the Theorem of Poincar\'e, to describe the dynamics of system (\ref{sistema}).
\medskip

\subsection{The characteristic polinomial  has no distinct roots with the same module}

Let $\lambda_1$, $\lambda_2$, $\lambda_3$ be the three roots of the characteristic polynomial of the linear difference
equation (\ref{ecuacionlineal}) in this case. A condition on the coefficients for this
case can be given by:
$$
\left( \displaystyle \frac{\frac{2}{3}\beta _{1}\alpha _{2}-\beta _{2}\alpha _{1}-\frac{2}{27}{\beta_{1}}^{3}\allowbreak }{2}\right)
^{2}\leq \left( \displaystyle \frac{\alpha _{2}+\frac{1}{3}{\beta_{1}}^{2}}{3}\right) ^{3},$$
with $\alpha_1 \ne 0$ or $\alpha_2 \leq 0$. Recall that we assume here that $\beta_2 \alpha_1 \ne \beta_1 \alpha_2$ and
$\beta_2\ne 0$.

If $\lambda_1$ is the characteristic root of maximal modulus, we will denote by $L$ the line
$$L=\{\ (x,y)\, :\,  \beta_2 x =
(\beta_1-\lambda_1)(y+\lambda_1) \ \}.$$

\begin{proposition} \label{distmod} Suppose that $\beta_2 \neq 0$ and every root of the characteristic polinomial of the linear difference
equation (\ref{ecuacionlineal}) is real and no two distinct roots have the same module. When $(x_{0},y_{0})$ is not in the forbidden
set, we have:
\begin{enumerate}
    \item If $|\lambda_1|>|\lambda_2|>|\lambda_3|$, then
    \begin{enumerate}
        \item the system (\ref{sistema}) admits exactly the three equilibria $\left(\displaystyle\frac{{\lambda_i}^2-\alpha_2}{\beta_2},\lambda_i\right)$, $i=1,2,3$;
        \item the fixed point $\left(\displaystyle\frac{{\lambda_1}^2-\alpha_2}{\beta_2},\lambda_1\right)$  attracts every complete solution starting on a point $(x_0,y_0)$ which does not belong to the line $L$;

        \item the corresponding solution to the system with initial condition $(x_0,y_0)\neq \left(\displaystyle\frac{{\lambda_3}^2-\alpha_2}{\beta_2},\lambda_3\right)$ and $(x_0,y_0)\in L$  converges to $\left(\displaystyle\frac{{\lambda_2}^2-\alpha_2}{\beta_2},\lambda_2\right).$
    \end{enumerate}
    \item If $|\lambda_1|>|\lambda_2|$ and $\lambda_1$ has algebraic multiplicity 2, then
    \begin{enumerate}
        \item the system (\ref{sistema}) admits exactly the two equilibria
        $\left(\displaystyle\frac{{\lambda_i}^2-\alpha_2}{\beta_2},\lambda_i\right)$, $i=1,2$;
        \item the fixed point $\left(\displaystyle\frac{{\lambda_1}^2-\alpha_2}{\beta_2},\lambda_1\right)$ attracts every complete solution except the other fixed point;
    \end{enumerate}
    \item If $|\lambda_1|>|\lambda_2|$ and $\lambda_2$ has algebraic multiplicity 2, then
    \begin{enumerate}
        \item the system (\ref{sistema}) admits exactly the two equilibria
        $\left(\displaystyle\frac{{\lambda_i}^2-\alpha_2}{\beta_2},\lambda_i\right)$, $i=1,2$;
        \item the fixed point $\left(\displaystyle\frac{{\lambda_1}^2-\alpha_2}{\beta_2},\lambda_1\right)$  attracts every complete solution starting on a point $(x_0,y_0)$ which does not belong to the line
        $L$;
        \item the corresponding solution to the system with initial condition $(x_0,y_0)\in L$ converges to
        $\left(\displaystyle\frac{{\lambda_2}^2-\alpha_2}{\beta_2},\lambda_2\right)$.
    \end{enumerate}

    \item If $\lambda_1$ has multiplicity 3, then
    \begin{enumerate}
        \item the system (\ref{sistema}) has a unique equilibrium $\left(\displaystyle\frac{{\lambda_1}^2-\alpha_2}{\beta_2},\lambda_1\right)$;
        \item the equilibrium is a global attractor.
    \end{enumerate}
\end{enumerate}
\end{proposition}
\dem  In all the cases, the equilibrium points are directly given by Proposition  \ref{sol2p}. The assertions concerning the asymptotic behaviour  can be derived as a consequence of Case 1 in  \cite[pg. 240]{elaydi}, bearing in mind that $$
x_n= \displaystyle \frac{v_{n+1}}{v_{n-1}}\displaystyle \frac{1}{\beta_2}-\displaystyle \frac{\alpha_2}{\beta_2}\, ,\quad
y_{n}=\displaystyle \frac{v_{n}}{v_{n-1}},
$$
and that  $v_n$ is the solution to the linear equation (\ref{ecuacionlineal}) with initial conditions $v_{-1}=1$, $v_0=y_0$ and $v_1=\beta_2 x_0 +\alpha_2$.$\square$

\subsection{The characteristic polinomial  has two distinct real roots with the same module}

It is easy to check that this case occurs when $\beta_1\ne 0$, $\beta_2\ne 0$, $\alpha_1=0$ and $\alpha_2>0$. Thus, the
roots of the characteristic polynomial of the linear difference equation (\ref{ecuacionlineal}) are $\beta_1$ and $\pm
\sqrt{\alpha_2}$.

\begin{proposition} \label{2real=mod} Suppose $\beta_1\ne 0$, $\beta_2\ne 0$, $\alpha_1=0$ and $\alpha_2>0$. Assume also that
$(x_{0},y_{0})$ is not in the forbidden set.
\begin{enumerate}
\item If ${\beta_1}^2= \alpha_2$, then

\begin{enumerate}
\item there are two equilibrium points $(0,\pm \beta_1)$;
\item the equilibrium point $(0,\beta_1)$ attracts every complete solution not starting on a point of the line $x=0$;
\item the solutions starting on a point $(x_0,y_0)$ of the line $x=0$ are prime period-two solutions except the two
equilibrium points $(0, \pm \beta_1)$.
\end{enumerate}

\item If ${\beta_1}^2> \alpha_2$, then
\begin{enumerate}
\item there are three equilibrium points $\left(\displaystyle \frac{{\beta_1}^2-\alpha_2}{\beta_2},\beta_1\right)$ and $(0, \pm
\sqrt{\alpha_2})$;
\item the equilibrium point $\left(\displaystyle\frac{{\beta_1}^2-\alpha_2}{\beta_2},\beta_1\right)$ attracts every
complete solution not starting on a point of the line $x=0$;
\item the solutions starting on a point $(x_0,y_0)$ of the line $x=0$ are prime period-two solutions except the two
equilibrium points $(0, \pm \sqrt{\alpha_2})$;
\end{enumerate}

\item If ${\beta_1}^2< \alpha_2$, then
\begin{enumerate}
\item there are three equilibrium points $\left(\displaystyle \frac{{\beta_1}^2-\alpha_2}{\beta_2},\beta_1\right)$ and $(0, \pm
\sqrt{\alpha_2})$;
\item the solutions starting on a point of the line $x=0$ are prime period-two solutions except the two
equilibrium points $(0, \pm \sqrt{\alpha_2})$;
\item the solutions starting on a point of the lines $\beta_2 x +
\displaystyle \frac{\alpha_2-{\beta_1}^2}{\beta_1}y=0$ or $x=\displaystyle \frac{{\beta_1}^2-\alpha_2}{\beta_2}$ are unbounded with the only exception of the fixed point $\left(\displaystyle \frac{{\beta_1}^2-\alpha_2}{\beta_2},\beta_1\right)$;
\item the solutions starting on any other point $(x_0,y_0)$ are bounded and each tends to one of the two-periodic solutions;
\end{enumerate}
\end{enumerate}
\end{proposition}
\dem   In all cases the affirmation \emph{a)} is a consequence of Proposition \ref{sol2p}.

When ${\beta_1}^2= \alpha_2$, the roots are $\beta_1$, with algebraic multiplicity two, and $-\beta_1$. By Proposition \ref{solb2no0} we know that any solution of the system can be written as
\begin{eqnarray*}
&\beta_2 x_n=\displaystyle \frac{(n+1) P_1 + P_2 + P_3 (-1)^{n+1}}{(n-1) P_1 + P_2 + P_3 (-1)^{n-1}} {\beta_1}^2 - {\beta_1}^2,&
\\
&y_n=\displaystyle \frac{n P_1 + P_2 + P_3 (-1)^n}{(n-1) P_1 + P_2 + P_3 (-1)^{n-1}} \beta_1,&
\end{eqnarray*}
where $P_1$, $P_2$ and $P_3$ actually satisfy 
\begin{equation}
\label{lsmev}
\begin{array}{rll}
 P_1 + P_2 - P_3  =   \displaystyle\frac{\beta_2 x_0 + {\beta_1}^2}{\beta_1}\, ,\quad
 P_2 + P_3 =   y_0\, ,\quad
-P_1 + P_2 - P_3   =  \beta_1.
\end{array}
\end{equation}

If $P_1\ne 0$, then $(x_n,y_n)$ obviously tends to $(0, \beta_1)$. From (\ref{lsmev}) we see that  $P_1=0$ if and only if  $x_0=0$ and, in such case,
 $x_n=0$ and $y_n$ takes alternatively the values $A \beta_1$ and $A^{-1} \beta_1$ with
$A=\frac{P_2+P_3}{P_2-P_3}$. Notice that $y_0\ne 0$ guaranties $P_2+P_3\ne 0$ and, since $\beta_1\ne 0$, we can not have
$P_1=0$ and $P_2-P_3=0$. This completes the proof of (1).

In the case $\beta_1^2\ne \alpha_2$, by Proposition \ref{solb2no0} we can write the general solution of the system as
\begin{eqnarray*}
&\beta_2 x_n=\displaystyle \frac{P_1 + \left [P_2 + P_3 (-1)^{n+1}\right ]\left ( \frac{\sqrt{\alpha_2}}{\beta_1}\right )^{n+1}}{P_1 +
\left [P_2 + P_3 (-1)^{n-1}\right ]\left ( \frac{\sqrt{\alpha_2}}{\beta_1}\right )^{n-1}} {\beta_1}^2 - \alpha_2,
&\\&
y_n=\displaystyle \frac{P_1 + \left [P_2 + P_3 (-1)^{n}\right ]\left ( \frac{\sqrt{\alpha_2}}{\beta_1}\right )^{n}}{P_1 + \left [P_2 +
P_3 (-1)^{n-1}\right ]\left ( \frac{\sqrt{\alpha_2}}{\beta_1}\right )^{n-1}} \beta_1,
&
\end{eqnarray*}
where $P_1$, $P_2$ and $P_3$ satisfy
\begin{equation}
\label{lsmev1}
\begin{array}{rll}
 P_1 \beta_1 + (P_2 - P_3)\sqrt{\alpha_2} & = &  \beta_2 x_0 + \alpha_2,\\
P_1+ P_2 + P_3 & =  & y_0,\\
P_1 {\beta_1}^{-1} + (P_2 - P_3)\sqrt{{\alpha_2}^{-1}} & = &  1.
\end{array}
\end{equation}

When $\beta_1^2>\alpha_2$ one immediately gets the results of statement (2) with an argument similar to that of the previous case. Therefore we will focus our
attention in the case $\beta_1^2<\alpha_2$. The condition $x_0=0$ is, according to (\ref{lsmev1}), equivalent to $P_1=0$ and in such case one gets $x_n=0$ and $y_n$ takes alternatively the values $K \sqrt{\alpha_2}$ and $K^{-1}
\sqrt{\alpha_2}$ with $K=\frac{P_2+P_3}{P_2-P_3}=\frac{y_0}{\alpha_2}$.  Now,
 if $P_1\ne 0$ and the initial conditions are taken such that $P_2+P_3\neq 0\neq P_2-P_3$, then $(x_n,y_n)$ tends obviously to  the 2-cycle $\{(0,K \sqrt{\alpha_2}),(0,K^{-1} \sqrt{\alpha_2})\}$ where $K=\frac{P_2+P_3}{P_2-P_3}$. On the contrary, if either $P_2+P_3= 0$ or $P_2-P_3=0$ (and only one of both equalities holds) then both sequences $x_n$ and $y_n$ are unbounded. From the system 
(\ref{lsmev1}) one gets that $P_2-P_3=0$ if and only if $x_0=(\beta_1^2-\alpha_2)/\beta_2$ and that $P_2+P_3= 0$ is equivalent to $\beta_2 x_0 +
\displaystyle \frac{\alpha_2-{\beta_1}^2}{\beta_1}y_0=0$. This shows the validity of {\it c)}.$\square$

\subsection{The characteristic
polinomial  has complex roots}

Now we consider the case in which the characteristic polynomial of the linear difference
equation has a couple of complex roots $\rho e^{\pm i\theta}$, with $\sin\theta>0$. Let $\lambda \neq 0$ be the real root.
It can be easily shown that
\begin{equation}\label{coef}
\beta_{1}=\lambda+2\rho\cos\theta; \quad \quad\alpha_{2}=-\left(  2\lambda
\rho\cos\theta+\rho^{2}\right) ; \quad \quad\beta_{2}\alpha_{1}=\lambda\rho
^{2}+\beta_{1}\alpha_{2}.
\end{equation}
and that this situation occurs when
$$
\left( \displaystyle \frac{\frac{2}{3}\beta _{1}\alpha _{2}-\beta _{2}\alpha _{1}-\frac{2}{27}{\beta_{1}}^{3}\allowbreak }{2}\right)
^{2} > \left( \displaystyle \frac{\alpha _{2}+\frac{1}{3}{\beta_{1}}^{2}}{3}\right) ^{3}.$$

By Proposition \ref{sol2p} we know that the unique equilibrium is
$\left(\displaystyle\frac{\lambda^2-\alpha_2}{\beta_2},\lambda\right)$. Denote by $L$ the line $$L=\{\ (x,y)\, :\,  \beta_2 x =
(\beta_1-\lambda)(y+\lambda) \ \}.$$ Notice that $(\beta_1-\lambda)(y+\lambda)=2y\rho\cos\theta-\alpha_2-\rho^{2}$. Also
observe that the equilibrium does not belong to $L.$
\begin{theorem} \label{complex}
Suppose $\beta_2\neq 0$ and the characteristic polynomial of the linear difference
equation to have complex roots and assume that $(x_{0},y_{0})$ is not in the forbidden set.
\begin{enumerate}
\item The solutions starting on the line $L$ remain on it and they are either all periodic or all unbounded.
\item If $|\lambda| > \rho$, then the unique equilibrium  attracts all the solutions not starting on $L$.
\item  If $|\lambda| < \rho$, then every non fixed bounded subsequence of a solution accumulates on $L$.
\item If $|\lambda| = \rho$, then every complete solution (neither starting on the fixed point nor on $L$) lies on a not degenerate conic, which does not contain the equilibrium.  \end{enumerate}
\end{theorem}
\dem  Assume that $(x_0, y_0)$ is not the fixed point. Using Proposition \ref{solb2no0}, we have%
\begin{eqnarray*}
\alpha_{2}+\beta_{2}x_{n} =\displaystyle \frac{P\lambda^{n+1}+2\rho^{n+1}\cos\left( a+\left( n+1\right)
\theta\right)  }{P\lambda^{n-1}+2\rho^{n-1}\cos\left( a+\left(  n-1\right) \theta\right)}\, ,
\\
y_{n} =\displaystyle \frac{P\lambda^{n}+2\rho^{n}\cos\left( a+n\theta\right) }{P\lambda^{n-1}+2\rho^{n-1}
\cos\left(  a+\left(  n-1\right)  \theta\right)  },
\end{eqnarray*}
where the constants $P\in\mathbb{R}$ and $a\in\lbrack0,2\pi)$, together with
$k\in\mathbb{R}^{+}$, are given by
\begin{equation}
\label{cofcomplex}
\left(
\begin{array}
[c]{ccc}%
\lambda & \rho e^{i\theta} & \rho e^{-i\theta}\\
1 & 1 & 1\\
1/\lambda &  e^{-i\theta}/\rho &  e^{i\theta}/\rho
\end{array}
\right)  \left(
\begin{array}
[c]{c}%
kP\\
ke^{ia}\\
ke^{-ia}%
\end{array}
\right)  =\left(
\begin{array}
[c]{c}%
\alpha_{2}+\beta_{2}x_{0}\\
y_{0}\\
1
\end{array}
\right).
\end{equation}

Observe that
we may consider $P\geq 0$, by replacing, if necessary, $a$ with $a +\pi$.

Let us consider the sequences $$\sigma_{n}=2\left(  \frac{\rho}{\lambda}\right) ^{n}\cos\left(
a+n\theta\right)  ; \quad \quad \quad \quad  \tau_{n}=2\left(  \frac{\rho}{\lambda}\right) ^{n}\sin\left(
a+n\theta\right).$$
It can be easily proved that
\begin{eqnarray}
\alpha_{2}+\beta_{2}x_{n}=\lambda^{2}\displaystyle \frac{P+\sigma_{n+1}}{P+\sigma_{n-1}}  , \quad \quad &
y_{n}=\lambda\displaystyle \frac{P+\sigma_{n}}{P+\sigma_{n-1}}\,, \label{solsigma}\\
\lambda\sigma_{n+1} =\rho\sigma_{n}\cos\theta-\rho\tau_{n}\sin\theta  , \quad \quad &
\rho\sigma_{n-1}=\lambda\sigma_{n}\cos\theta+\lambda\tau_{n} \sin\theta\, .\label{sigmatau}
\end{eqnarray}
As a consequence, $\lambda^2 \sigma_{n+1}-2\lambda \rho \sigma_n \cos \theta+\rho^2\sigma_{n-1}=0$, and then
\[
\alpha_{2}+\beta_{2}x_{n}=2\rho y_{n}\cos\theta-\rho^{2}+P\displaystyle \frac{\lambda
^{2}-2\rho\lambda\cos\theta+\rho^{2}}{P+\sigma_{n-1}}\, ,
\]
which is equivalent to
\begin{eqnarray} \label{casirecta}
\beta_{2}x_{n}-(\beta_1-\lambda)(y_n+\lambda)=P\displaystyle \frac{\lambda
^{2}-2\rho\lambda\cos\theta+\rho^{2}}{P+\sigma_{n-1}}.
\end{eqnarray}
Using (\ref{cofcomplex}) one has that $(x_0,y_0)\in L$ if and only if $P=0$ and from (\ref{casirecta}) we then get that
 $(x_n, y_n)\in L $ for all $n\geq 1$.

Furthermore, by (\ref{solsigma}), we see that if $(x_0,y_0)\in L$, then the solution $(x_n, y_n)$ is periodic whenever
${\theta}/{\pi}$ is a rational number and unbounded otherwise.

Assume now that the solution $(x_n, y_n)$ does not start on $L$, this to say, $P\neq 0$. We will now distinguish the three
cases: $|\lambda| > \rho$, $|\lambda| < \rho$ and $|\lambda| = \rho$.

If $|\lambda| > \rho$, then by (\ref{solsigma}) one immediately has $x_n\to \frac{\lambda^2-\alpha_2}{\beta_2}$ and $y_n \to
\lambda$.

Suppose now that $|\lambda| < \rho$. If $(x_{n_k}, y_{n_k})$ is a subsequence satisfying that $\inf\limits_{k}\left|  \cos\left(  a+\left(
n_{k}-1\right)  \theta\right)  \right|  >0,$ then one obviously has $\sigma_{n_{k}-1}\rightarrow\infty$. Using the definition of $\sigma_n$, one easily gets that 
$\displaystyle \frac{\sigma_{n_{k}}}{\sigma_{n_{k}-1}}$ is bounded. Then, $\left(  x_{n_{k}},y_{n_{k}}\right) $ is a
bounded subsequence and equation (\ref{casirecta}) shows that it is attracted by the line $L$.

On the other hand, if $\cos\left(
a+\left(  n_{k}-1\right)  \theta\right)  \rightarrow 0,$ then the left equation in (\ref{sigmatau}) lead us to $\left|
\sigma_{n_{k}}\left(  \frac{\lambda}{\rho}\right)  ^{n_{k}}\right|
\rightarrow2\sin\theta>0.$ Thus,  $\sigma_{n_{k}}\rightarrow\infty$ and using (\ref{sigmatau}) once more we get
$\displaystyle \frac{\sigma_{n_{k}}}{\sigma_{n_{k}-1}}\rightarrow\infty$. Therefore,
$\left(  x_{n_{k}},y_{n_{k}}\right)  $ is an unbounded subsequence.

Finally, let us suppose $\rho=|\lambda|$. If we consider the change of variables
\begin{eqnarray*}
\overline{x}  &  =&\left(  \beta_{2}x+\alpha_{2}-\rho^{2}\right)
\displaystyle\frac{\lambda}{2\lambda\cos\theta-2\rho}-\left(  y-\lambda\right)  \displaystyle\frac{\rho\lambda\cos\theta}
{\lambda\cos\theta-\rho}\, ,\\
\overline{y}  &  =&\left(  \beta_{2}x+\alpha_{2}-\rho^{2}\right)
\displaystyle\frac{1}{2\sin\theta}-(y-\lambda)\displaystyle\frac{(\lambda+\rho\cos\theta)}{\sin\theta}
\,,
\end{eqnarray*}
then one may deduce from (\ref{sigmatau}) that $\overline{x_n}=\displaystyle \frac{\rho \lambda\sigma_{n-1}}{P+\sigma_{n-1}},
\overline{y_n}=\displaystyle \frac{ \rho \lambda\tau_{n-1}}{P+\sigma_{n-1}}$. Therefore, one immediately gets that
$$\overline{x_n}^{2}+\overline
{y_n}^{2}=4\frac{(\rho\lambda)^2}{(P+\sigma_{n-1})^2}\, ,\quad (\overline{x_n}-\rho\lambda)^2=P^2\frac{(\rho\lambda)^2}{(P+\sigma_{n-1})^2}\, ,$$
which clearly shows that $(\overline{x_n},\overline{y_n})$ lies in the conic $\overline{x}^2+\overline{y}^2=\frac{4}{P^2}(\overline{x}-\rho\lambda)^2$, having its focus  in $(0,0)$, its directrix in the line $\overline{x}=\rho\lambda$ and eccentricity  $2/P$.  
Further, one immediately sees that the fixed point $\left(\displaystyle\frac{\lambda^2-\alpha_2}{\beta_2},\lambda\right)$ is transformed by the change of variables above in $(0,0)$ and, hence, it does not belong to the conic. $\square$

\begin{remark}
In the case $|\lambda|<\rho$ of this last theorem, one might conjecture that every subsequence  of a solution (even a non-bounded one) actually approaches the line $L$; but this is not the case. Let us take, for example, the system whith $\alpha_1=1,\beta_1=3,\alpha_2=-4,\beta_2=-10$, in which the characteristic roots of the associated polynomial are given by $\lambda =1$ and $\sqrt{2}e^{i\pi/4}$ and consider the solution starting on $(x_0,y_0)=(-11/20,3/2)$. We then have that $a=0$, $P=1$ and $\sigma_{2+4k}=0$ for all $k\ge 0$. One may use equation (\ref{casirecta}) to show that all the points of the form $(x_{3+4k},y_{3+4k})$ lay on the line $10x+2y+3=0$ while the line $L$ is given by $10x+2y+2=0$. Note, however, that the subsequences $(x_{4k},y_{4k}),(x_{1+4k},y_{1+4k})$ and $(x_{2+4k},y_{2+4k})$ are all bounded and converge respectively to $(-3/5,2),(-2/5,1)$ and $(-1/5,0)$, which do belong to $L$.

It should also be noticed that the fixed point lays on the line $10x+2y+3=0$. This is also the case in the general setting. It follows from (\ref{casirecta}) that whenever $\sigma_{n_k-1}=0$ then the point $(x_{n-k},y_{n-k})$ is on the line containing the fixed point which is parallel to $L$. 
\end{remark}

\begin{remark} Notice that, according to the results in \cite{b-l}, when $|\lambda|=\rho$ and the argument $\theta$ of the complex root is a rational multiple of $\pi$, the system is globally periodic.
\end{remark}

\subsection{Stability of fixed points}

We finish this section with the complete study of the stability of the fixed points in the case $\beta_2\neq 0$. 

\begin{theorem} Suppose that $\beta_2\neq 0$, let $\lambda$ be a real eigenvalue of the matrix $A$ given in (\ref{matriz}). Let $\left(\frac{\lambda^2-\alpha_2}{\beta_2},\lambda\right)$ be the associated fixed point and denote by $\rho (A)$ the spectral radius of $A$.
\begin{enumerate}
\item If $|\lambda|<\rho (A)$, then the associated fixed point is unstable.
\item If $|\lambda|=\rho (A)$, then the associated equilibrium is stable if and only if every eigenvalue whose modulus is $\rho (A)$ is a simple eigenvalue. Moreover, the stability is asymptotic if and only if $\lambda$ is a simple eigenvalue and it is the unique eigenvalue of $A$ whose modulus is $\rho (A)$.
\end{enumerate}\end{theorem}
\dem  The first statement was already proved in Proposition \ref{esta-linear}. 
Besides, in such Proposition, we have shown that if every eigenvalue whose modulus is $\rho (A)$ is simple then the associated equilibrium is stable. Let us prove the converse. 

According to the results of the previous subsections, the only cases in which one has a non-simple eigenvalue of maximal modulus are the cases treated in Proposition  \ref{distmod} (1) and (4) and the first case of Proposition \ref{2real=mod}. We will see that in such cases the equilibrium points associated to eigenvalues of maximal modulus are unstable. 

We begin with the case of an eigenvalue $\lambda_1$ of maximal modulus with multiplicity 2. For each  $N\in {\mathbb N}$, $N>1$ one may consider the solution with initial conditions $(x_0,y_0)=(\frac{\lambda_1^2-\alpha_2}{\beta_2}-\frac{2\lambda_1^2N}{(N^2+1)\beta_2},\lambda_1-\frac{\lambda_1N}{N^2+1})$. The solution of  (\ref{ecuacionlineal})  in such case is given by $v_n=\lambda_1^{n+1}(N^2+1-Nn-N)/(N^2+1)$, which cannot vanish since $N>1$.
For this solution one has $|y_{N}-\lambda_1|=|\lambda_1|N$, proving that the equilibrium $(
\frac{\lambda_1^2-\alpha_2}{\beta_2},\lambda_1)$ is unstable.

Similarly, if $A$ has a unique eigenvalue $\lambda$ of multiplicity 3 then, for each $N\in{\mathbb N}$, $N\neq 0$ let us consider $(x_0,y_0)=(
\frac{\lambda^2-\alpha_2}{\beta_2}-2\frac{\lambda^2}{\beta_2 N^2},\lambda)$. The corresponding solution to (\ref{ecuacionlineal}) is given by
$v_n=\frac{(N^2-n-n^2)\lambda^{n+1}}{N^2}$. It is not difficult to see that $v_n\neq 0$ for all $n\geq 1$ and then the solution to our system (\ref{sistema}) is complete. Further, since $y_n=v_n/v_{n-1}$ one gets that $|y_N-\lambda|=2|\lambda|$. Therefore, the fixed point $(
\frac{\lambda^2-\alpha_2}{\beta_2},\lambda)$ is not stable.

When $\alpha_1=0,\beta_1^2=\alpha_2\neq 0$, there are two equilibrium points associated to eigenvalues of maximal modulus: $(0,\pm\beta_1)$. The fixed point $(0,-\beta_1)$ is, according to the result of Proposition \ref{2real=mod}, unstable since the other equilibrium attracts all the solutions not starting on the line $x=0$. To see that $(0,\beta_1)$ is also unstable,  let us choose, for each odd $N\in{\mathbb N}$, the solution starting at $(x_0,y_0)=(\frac{-2\beta_1^2}{N\beta_2},\beta_1)$. Then, using equation (\ref{lsmev}) and the expression for $y_n$ given just above such equation, we have 
$v_n=\beta_1^{n+1}+P_1n\beta_1^n$ if $n$ is even and $v_n=\beta_1^{n+1}+P_1(n+1)\beta_1^n$ if $n$ is odd,
where $P_1=-\beta_1/N$. Since $N$ is odd, we see that $y_n$ exists for all $n\in {\mathbb N}$ and, further, we get that
$|y_{N}-\beta_1|=2|\beta_1|$, which clearly implies that $(0,\beta_1)$ cannot be stable.

Finally, it only remains to prove that when $A$ has distinct simple eigenvalues whose modulus equal $\rho (A)$ then the fixed point is not asymptotically stable. But this situation can only happen if either one has the situation described in Proposition \ref{complex} (4) or the one given in Proposition \ref{2real=mod} (3). In the case of complex eigenvalues we had seen that all the orbits lie on  conics not going through the fixed point and, hence, it cannot be asymptotically stable. In the other case, it is clear that the fixed points $(0,\pm\sqrt{\alpha_2})$ are not attracting since every solution starting on the line $x=0$ is 2-periodic.
$\square$

\begin{remark} It is interesting to notice that in the three cases in which there is an eigenvalue of maximal modulus with multiplicity larger than 1, the corresponding fixed point is attracting but unstable.
\end{remark}

\section{Non-negative solutions to the system with non-negative coefficients}

When the coefficients of our system (\ref{sistema}) are non-negative and we restrict ourselves to non-negative initial conditions, many of the cases studied in the previous sections cannot appear. Further, in such case one may describe which kind of orbits appear and their assymptotic behaviour without the previous calculation of the characteristic roots.

It shoud be noticed that whenever the coefficients in system  (\ref{sistema}) are non-negative and $\alpha_1\beta_2\neq\alpha_2\beta_1$, every initial condition $(x_0,y_0)$ with $x_0\ge 0$, $y_0> 0$ gives rise to a complete orbit except for $\alpha_2=0$ where the condition $x_0> 0$ is also necessary.

It will be convenient to study independently the case $\alpha_1\beta_2=0$. Next result is a simple summary of the results in Section 3 and Proposition \ref{2real=mod} and, hence, we omit its proof.

\begin{corollary} Let us consider that the coefficients in system  (\ref{sistema}) are non-negative and $\alpha_1\beta_2=0\neq\alpha_2\beta_1$.
\begin{enumerate}
\item If $\beta_2=0$, one has:
\begin{enumerate}
\item When $\alpha_2\leq\beta_1^2$ there are no non-negative periodic orbits and all non-negative solutions are unbounded, with the only exception of the  case
$\alpha_2=\beta_1^2$, $\alpha_1=0$, which is globally 2-periodic.
\item When $\alpha_2>\beta_1^2$ there exists a non attractive fixed point $\left(\frac{\alpha_1}{\sqrt{\alpha_2}-\beta_1},\sqrt{\alpha_2}\right)$ and the whole line $(\alpha_2-\beta_1^2)x_0=\alpha_1(\beta_1+y_0)$ of 2-periodic solutions. Every other non-negative solution is bounded and converges to one of the 2-cycles.
\end{enumerate}
\item If $\beta_2\neq 0=\alpha_1$, then every non-negative solution is bounded and the ones starting in the line $x_0=0$ are 2-periodic. Morevoer,
\begin{enumerate}
\item When $\alpha_2<\beta_1^2$ there are two non-negative fixed points: $\left(\frac{\beta_1^2-\alpha_2}{\beta_2},\beta_1\right)$, which attracts all non-periodic non-negative solutions, and $(0,\sqrt{\alpha_2})$.
\item When $\alpha_2=\beta_1^2$ there is a unique non-negative equilibrium $(0,\beta_1)$ which attracts all non-periodic non-negative solutions.
\item When $\alpha_2>\beta_1^2$ the unique non-negative equilibrium is $(0,\sqrt{\alpha_2})$ which is not an attractor. Every non-negative solution converges to one of the periodic solutions.
\end{enumerate}
\end{enumerate}
\end{corollary}

The remaining cases are jointly treated in the following result. All the definitions and results on non-negative matrices which are used in its proof may be found in \cite[Ch. 8]{hj}

\begin{proposition} Suppose that system  (\ref{sistema}) has  non-negative coefficients and that $\alpha_1\beta_2\neq 0$.
\begin{enumerate}
\item If $\alpha_2\neq 0$ or $\beta_1\neq 0$ then there is a unique non-negative (actually, positive) stable equilibrium which attracts all non-negative solutions.
\item If $\alpha_2=\beta_1=0$, the system is globally 3-periodic with a unique equilibrium.
\end{enumerate}
\end{proposition} 
\dem 
Let us consider $A$ as in (\ref{matriz}). A simple calculation shows that $(A+I)^2$ is positive and, therefore, $A$ is irreducible. Then the spectral radius $\rho (A)$ is a strictly positive simple eigenvalue of $A$. 

 If there exists another eigenvalue $\lambda$ such that $|\lambda|=\rho (A)$ then, since $A$ is non-negative and irreducible, the eigenvalues of $A$ should be $\lambda_{k+1}=\rho (A)e^{ik\pi/3}$ where $k=0,1,2$ and, consequently, $A^3=\rho (A)^3I$. The direct computation of $A^3$ shows that this is possible if and only if  $\alpha_2=\beta_1=0$ and, hence, in that case, the system is 3-periodic and the only equilibrium is the one associated to the real eigenvalue $\rho (A)$. 

In the remaining cases, $\lambda_1=\rho (A)$ is a dominant eigenvalue and, according to our results of Proposition \ref{esta-linear}, Proposition \ref{distmod} and Theorem \ref{complex}, the corresponding fixed point is stable and attracts all complete solutions except those starting on the line
$$L=\{\ (x,y)\, :\,  \beta_2 x =(\beta_1-\lambda_1)(y+\lambda_1) \ \}.$$
Since $\lambda_1$ is the largest eigenvalue of $A$, one has that $\mbox{det}(A-\mu I)<0$ for all $\mu>\lambda_1$. However, $\mbox{det}(A-\beta_1 I)=\alpha_1\beta_2>0$, showing that $\beta_1<\lambda_1$. Thus, for every $x_0\geq 0$ and $y_0> 0$ one obtains
$\beta_2x_0\ge 0$ and $(\beta_1-\lambda_1)(y_0+\lambda_1)< 0$, which proves that $(x_0,y_0)\not\in L$.

The equilibrium associated to the eigenvalue $\lambda_1=\rho (A)$ is $\left(\frac{\lambda_1^2-\alpha_2}{\beta_2},\lambda_1\right)$, which is positive since, as before, one sees that $\mbox{det}(A-\sqrt{\alpha_2} I)=\alpha_1\beta_2>0$ and hence $\lambda_1>\sqrt{\alpha_2}$. 
$\square$ 

\section*{Acknowledgements} 
We want to thank Professor Eduardo Liz for his useful comments and suggestions. This work was partially supported by  MEC Project MTM2007-60679.

\bigskip

\noindent {\bf Authors' address:}\\ 
I. Bajo, Depto. Matem\'atica Aplicada II, E.T.S.E.
Telecomunicaci\'on, Campus Marcosende, Universidade de Vigo, 36310 Vigo, Spain. ibajo@dma.uvigo.es\\
D. Franco and J. Per\'an, Departamento de
Matem\'atica Aplicada, E.T.S.I. Industriales, UNED, c/ Juan del Rosal 12, 28040 Madrid,
Spain. dfranco@ind.uned.es, jperan@ind.uned.es
\end{document}